\documentclass[a4paper,12pt]{article}
\usepackage{amsmath,amsfonts}
\usepackage[english]{babel}
\usepackage{a4wide}

\usepackage{amssymb,amsthm}

\usepackage{color}

\usepackage[cp1251]{inputenc}

\theoremstyle{plain}

\newtheorem{thm}{Theorem}
\newtheorem{lem}{Lemma}
\newtheorem{cor}{Corollary}
\theoremstyle{remark}
\newtheorem{rem}{Remark}
\theoremstyle{definition}
\newtheorem{Stat}{Statement } 

\newcommand{\E}{\mathbb{E}}

\title{Estimates of the probability of a regenerative process reaching a high level}

\author{Kateryna~Akbash$^a$, Ivan~Matsak$^b$, Oleg~Zakusylo$^b$}
\date{}

\begin{document}
\maketitle

\begin{center}
$^a$ Volodymyr Vynnychenko Central Ukrainian State
University,\\
$1$, Shevchenko street, Kropyvnytskyi 25006, Ukraine\\
\medskip
$^b$ Taras Shevchenko National University of Kyiv,\\
$2/6 $, Academician Glushkov avenue, Kyiv 03127, Ukraine\\
\medskip
kateryna.akbash@gmail.com (K.S.Akbash),  zakusylo@knu.ua
(O.K.~Zakusylo), i.m.k@ukr.net (I.K.Matsak), 
\end{center} \bigskip

\begin{abstract}
The problem of estimating the probability of a random process reaching a certain level is well known. In this article, two-sided estimates are established for the probability that a regenerative process reaches a high level. Two auxiliary results for geometric sums with delay will play an important role. Examples of application to random processes describing queue lengths in queueing theory are also given.
\end{abstract}

 \medskip

\bigskip

\section{Introduction. Main results}
The problem of estimating the probability of a random process (r.p.) reaching a certain level arises naturally in applications. For example, in mathematical reliability theory, the problem of estimating the distribution of the time to failure of a redundant system with recovery is well known \cite{gbs}. Similar problems also exist in queueing theory, inventory management, etc. Such problems have been intensively studied in (\cite{ADM_24}, \cite{asm_2}, \cite{gl_kou}, \cite{gbs},
\cite{Ka97}, \cite{v_ka}, \cite{OK_IK_20}).

Clearly, this topic is closely related to the study of extreme values of random processes. Moreover, the extreme values of processes that describe queue lengths or queue wait times for queueing systems have also been studied in many studies \cite{asm_2}, \cite{coh}, \cite{igl}, \cite{OK_IK_20}.

It should be noted that in many cases, the corresponding mathematical models can be adequately described by regenerative random processes. In this article, we restrict our attention to regenerative processes. The solution to our problem will be based on limit theorems and corresponding estimates for geometric sums (see, for example, \cite{ADM_24}, \cite{bro90}, \cite{Ka97}, \cite{pe_r11}).

Let us recall the definition of a regenerative r.p. Let $T$ be a non-negative random variable (r.v.)  such that  $\mathbf{P}( T > 0)>0$,  $\mathbf{P}( T < \infty) =1 $, and let a r.p.  $\xi (t)$,
be given on the interval $[0,T),  \xi (t)\in \mathbb{R}^1 $. Such a pair $\pounds = (T, \xi (t))$, $t\in [0,T)$ is called a cycle. Note that the r.v.  $T$  and the r.p. $\xi (t)$ may be dependent.

Let $\pounds_k = (T_k , \xi_k (t))$,  $k\geq 1 $, be an infinite sequence of independent cycles, identically distributed with the cycle $\pounds = (T,\xi(t))$. Then the process $X(t), t\geq 0, $ is regenerative if it can be represented as follows:
\[
X(t) = \xi_k (t- S_{k-1}),  \quad  \mbox{for} \quad   t \in [S_{k-1}
, S_{k}) ,
\]
where $S_{k} = T_1 + \ldots + T_k , k \geq 1$, $S_{0}=0$,  (see, for example,  \cite{ks}, part II, ch.2, \cite{fe2}, ch.11, \S 8). Points $S_{k}$ are called moments of regeneration, and the interval $[S_{k-1} , S_{k})$ is called the $k$th period of regeneration. Assume that  $X(S_k) = \xi_k (0) =0, k\geq 0 $, that is,  the process  $X(t) $ at the beginning of the regeneration cycle starts from state $0$.

We impose a separability condition on the r.p.  $ \xi_i (t)$. Then we can introduce the r.p.
\[
\bar{X}(t) = \sup_{0\leq s < t} X(s), \quad t\geq 0
 \]
and consider a sequence of independent identically distributed random variables (i.i.d.r.v)

\[
 \bar{X}_k = \sup_{S_{k-1} \leq s <
S_k } X(s) , \quad k=1, 2, \ldots ,
\]
i.e., $\bar{X}_k $  is the maximum value of process ${X}(t)$ at the $k$th regeneration cycle.

Let's introduce some necessary notations. We divide the regeneration cycles into two types. The first type includes cycles in which some event does not occur, and the second type includes cycles in which event $A$ occurs. Event $A$ is interpreted as a failure of some system, which is described by process ${X}(t)$.

In applied problems, ${X}(t)$ often denotes the number of faulty elements of a system at time $t$. Then event $A$ occurs on the $k$th cycle if at some point  $t \in [S_{k-1},  S_k )$ all $u$ working elements of the system fail, i.e.
\begin{eqnarray}\label{f1.1}
A = A_u = \{\bar{X}_k \geq u \}.
 \end{eqnarray}
In this article, by $A_u$ we will mean an event that is defined by the formula (\ref{f1.1}), where $u>0$ is a fixed positive constant, $\bar{A_u}=\{\bar{X}_k
< u \}$ denotes the complementary event.

Let
\[
 \quad  q = q(u) = \mathbf{P}(A_u) =\mathbf{P}(\bar{X}_k \geq u) > 0
 .
 \]
Set  $ q^{*}=-\ln(1-q)$. It will be more convenient for us to normalize the corresponding amounts by the value $ q^{*}$, which is equivalent to $q$.
More precisely,
\[
q \leq  q^{*} \leq q+ q^2  \quad  \mbox{for} \quad  0< q < 1/2.
\]
It will be shown later that our estimates are meaningful only for small values of $q$. Therefore, we will assume that throughout $0< q < 1/2$.

Let us denote by $ \mathbb{F}(x)= \mathbf{P}( T < x)$ the distribution function of the r.v.  $T$, and by $T_k^{-} $ the lengths of the cycles on which event $\bar{A_u}$ occurred (cycle type $1$),
\[
\mathbb{F}^{-} (x)= \mathbf{P}( T_k^{-} < x) = \mathbf{P}( T_k < x /
\bar{A_u} ),
\]
\[
m_1^{-} = \int_0^{\infty} x d\mathbb{F}^{-} (x) , \quad m_{r}^{-} =
\int_0^{\infty} x^{r} d\mathbb{F}^{-} (x) .
\]
Define analogously
\[
T_k^{+},  \quad   \mathbb{F}^{+} (x) , \quad m_1^{+} , \quad
m_{r}^{+}
\]
to type-$2$ cycles, in which event $A_u$ occurred.
In addition, for cycles of the $2$nd type, we will consider the value $\hat{m}_1^{+}$ as the mean time for the process $\xi (t)$ to reach the level $u$  in the corresponding cycle.

Finally, we introduce the r.v. $\mathfrak{T}_X (u)$ - the first moment when process $X (t)$ reaches level $u$, that is,
$$\mathfrak{T}_X (u) = \inf(t\geq 0: \quad X (t) \geq u), $$
and also
\[
 \quad \Delta_X (x) =  1 - \exp(-
{x}) -   G_{X} \left(\frac{x m_1^{-}}{q^*}\right),  \quad G_{X} (x) =
  \mathbf{P}(  \mathfrak{T}_X (u) \leq x).
\]

We impose the following condition on the process $X(t)$: for some ${\gamma > 2}$
\begin{eqnarray}\label{f1.2}
m_{\gamma} = \mathbf{E}T^{\gamma} < \infty ,
 \end{eqnarray}
and thus $m_1 = \mathbf{E}T < \infty$ and  $m_2 = \mathbf{E}T^2 < \infty$.

The main result of this article is formulated in the following theorem.
\begin{thm}\label{t1.1} Let condition (\ref{f1.2}) be satisfied,
$0< x < 1$ is a fixed number and $q > 0$. Then

\begin{eqnarray}\label{f1.3}
   &&\exp(-x){q^*} \left(\frac{\hat{m}_1^{+}}{m_1^{-}} - \mathbb{C}_0 (\mathbb{F}, x)\right)
    \leq  \Delta_X (x)  \leq    \exp(-x){q^*} \times   \nonumber \\
   & \times &
   \left(\frac{\hat{m}_1^{+}}{m_1^{-}}(1+2x(e-2))
 +     \mathbb{C}_1 (\mathbb{F}, x) + \mathbb{C}_2 (\mathbb{F}, x) \right)
 .
  \end{eqnarray}
where
\[
 \mathbb{C}_0 (\mathbb{F}, x)= \frac{m_2^{-}}{(m_1^{-})^2} + \frac{\gamma m_{\gamma}(q^*)^{\gamma-2}(1+q)}{(\gamma-1)x^{\gamma-1}
 (m_1^{-})^{\gamma}},
\]
\[
 \mathbb{C}_1 (\mathbb{F}, x)= \left(2x(e-2) +1\right)\frac{m_{\gamma }(q^*)^{\gamma -1}(1+q)}{(m_1^{-})^{\gamma} x^{\gamma
 -1}}
 ,
\]
\[
 \mathbb{C}_2 (\mathbb{F},x) = (e-2)\left(\left(\frac{2m_2^{-} }{(m^{-}_1)^2}  +1\right) \left(x+\left(\frac{ m^{-}_2
  }{(m^{-}_1)^2}-1\right) q^{*}\right)+
  \frac{2m^{-}_2 x}{(m^{-}_1)^2 } -2x + q^{*}\right)  .
\]

\end{thm}

\begin{rem}\label{z1}
It can be assumed that in many practical problems the probability of event
$A$ (the probability of system failure) is quite small, and therefore the estimate
(\ref{f1.3}) is most interesting for small $x$. Therefore, we restrict ourselves to the case when $0 < x < 1$. It should be noted that the left-hand side inequality in
(\ref{f1.3}) will remain true for all $ x > 0$.

\end{rem}

\begin{cor}\label{n1.1}
If, under the conditions and notations of Theorem \ref{t1.1}  the distribution function
$\mathbb{F}(x)$ is fixed and does not depend on the parameter $u$,
 $q(u) \rightarrow 0$ when $u\rightarrow \infty$, then

\begin{eqnarray}\label{f1.4}
   &&\exp(-x){q^*} \left(\frac{\hat{m}_1^{+}}{m^{-}_1} - \frac{m_2}{m^2_1}+o(1)\right)
    \leq  \Delta_X (x)  \leq    \exp(-x){q^*} \times   \nonumber \\
   & \times &
   \left(\frac{\hat{m}_1^{+}}{m_1^{-}}(1+2x(e-2))
  +   (e-2)\left(\frac{4m_2 x }{m_1 ^2} -x +o(1) \right)\right).
  \end{eqnarray}
\end{cor}

\begin{rem}\label{z2}
Under the conditions of Corollary \ref{n1.1}  the lower bound for function $1-
G_{X} (\frac{x m_1^{-}}{q^*})$ can be written as follows: for
$u \rightarrow \infty$
\begin{eqnarray}\label{f1.5}
1- G_{X} \left(\frac{x m_1^{-}}{q^*}\right) \geq \exp(-x)  + \exp(-x){q^*}
\left(\frac{\hat{m}_1^{+}}{m^{-}_1} - \frac{m_2}{m^2_1}+o(1) \right).
\end{eqnarray}
This is a pessimistic estimate. In practice, it is often the pessimistic assessments that are often of primary importance (see \cite{gbs}).

It is clear that the estimate (\ref{f1.5})
is asymptotically more accurate than the estimate obtained in \cite{ADM_24}:\\
$$
1- G_{X}\left(\frac{x m_1^{-}}{q^*}\right) \geq  \exp(-x) -
  \exp(- {x})q^* \times
$$
$$
 \times \left((2x(e-2)+1)\frac{{m}^{+}_1  }{m^{-}_1} + \\
   (e-2)  \left(\frac{4 m_2 }{(m^{-}_1)^2} -1\right)  +o(1)\right).
$$

The same applies to some uniform estimates (see, for example, Theorem 4 of \cite{v_ka} and the simple and rather sharp estimate for geometric sums with delay in Theorem 2.2 of  \cite{bro90}).

Note also that under the conditions of Corollary \ref{n1.1} for
$\frac{\hat{m}_1^{+}}{m^{-}_1} > \frac{m_2}{m^2_1}$, for instance, if $\frac{\hat{m}_1^{+}}{m^{-}_1}  \rightarrow \infty$ (see example 2 in section 4), the approximation
\[
1- G_{X} \left(\frac{x m_1^{-}}{q^*}\right) \sim 1- \exp(-x)  + \exp(-x){q^*}
\frac{\hat{m}_1^{+}}{m^{-}_1}
  \]
is asymptotically more accurate than the traditional approximation:
\[
 1-G_{X} \left(\frac{x m_1^{-}}{q^*}\right) \sim  \exp(-x). \]
 It should be noted, however, that in the general case the problem of finding the value of $\hat{m}_1^{+}$ is rather involved.

\end{rem}

\begin{rem}\label{z3}
Consider a regenerative process of a special type, in which any regeneration period consists of two independent parts: $T= T^{'}+ T^{''}$,  where the first part  $T^{'}$ has an exponential distribution
\[
\mathbf{P}(T^{'}> t) = \exp(-\lambda t),
\]
and the second part $T^{''}$ has an arbitrary distribution,  $\mathbf{E}T^{''}< \infty$.
For such a special class of regenerative processes, the following simple pessimistic estimate is known \cite{gbs}: $\forall x>0$
\begin{eqnarray}\label{f1.6}
1- G_{X} \left(\frac{x }{\lambda q}\right) \geq \exp(-x) .
\end{eqnarray}

It is easy to see that for the general regenerative process under the conditions of Corollary \ref{n1.1} the following  inequality follows from estimate (\ref{f1.5}):
for small $q$
\[
1- G_{X} \left(\frac{x m_1^{-}}{q^*}\right) \geq  \exp(-x) . \]

Since $1/\lambda < m_1^{-}$, the the latter estimate has certain advantages over estimate (\ref{f1.6}).
\end{rem}

Thus, Theorem \ref{t1.1} and Corollary \ref{n1.1}
clarify the results of \cite{ADM_24} and related works.

In the following Section 2, two lemmas for geometric sums with delay are established, which will play an important role in the proof of Theorem \ref{t1.1}. The proof of Theorem \ref{t1.1} itself is given in Section 3. Finally, at the end of the article, Section 4 discusses some examples of using the obtained results to study processes from queueing theory.

\section{Two lemmas for geometric sums}

Let $\zeta , \zeta_1, \zeta_2, \ldots  $ be a sequence of positive i.i.d.r.v. with distribution function $F(x)=
\mathbf{P}(\zeta \leq x)$, $\nu $ be a geometrically distributed random variable independent of sequence $(\zeta_i )$,
\[
 \mathbf{P}( \nu = k) = q (1-q)^{k-1}, \quad  k=1, 2, \ldots , \quad
0<q<1 . \]
Let's put
\[
S_\nu = \sum_{k=1}^{\nu} \zeta_k .
\]
Such sums $S_\nu $ are called geometric random sums.
It is well known that under the condition $ a_1 =  \mathbf{E}\zeta <  \infty$  the following asymptotic relation holds \cite{Ka97}:\\

\begin{eqnarray}\label{f2.1}
 \forall x>0 \quad  \lim_{q \to 0} \mathbf{P}\left( \frac{q}{a_1} S_\nu \leq x\right)
= 1 - \exp(- {x}).
  \end{eqnarray}
It should be noted that geometric sums have been studied extensively. The main effort has been directed toward estimating the rate of convergence in equality (\ref{f2.1}) (see \cite{ADM_24}, \cite{bro90},
\cite{Ka97}, \cite{pe_r11},  \cite{v_ka}).

Thus, in  \cite{bro90} and  \cite{Ka97} under broad conditions on the r.v., fairly uniform estimates of the type
\begin{eqnarray}\label{f2.2}
\sup_{x\geq 0}  | \Delta_0 (x)| \leq C_0 q ,
  \end{eqnarray}
were established, where the constant $C_0$ depends on the distribution of the r.v. $\zeta$,
\[
 \quad \Delta_0 (x) =  1 - \exp(-
{x}) - \mathbf{P}\left( \frac{q}{a_1} S_\nu \leq x \right).
\]

Some continuous estimates are also known (see \cite{ADM_24},
 \cite{Ka97}, \cite{v_ka}).

A detailed overview of results for geometric sums and their applications, as well as a complete bibliography on this topic, can be found in the book  \cite{Ka97} and article \cite{pe_r11}.

For our purposes in what follows, the main case will be the case of geometric sums with delay  $S^{(d)}_\nu $, that is, when

$$ S^{(d)}_n = \zeta_1^{(d)} + \zeta_2 + \ldots\ldots + \zeta_n,  $$
where the first term $ \zeta_1^{(d)} $ has a distribution
$F^{(d)} (x) = \mathbf{P}(\zeta_1^{(d)} \leq x)$ different from $F(x)$.

Let us assume that r.v. $\zeta_1^{(d)},  \zeta_2, \zeta_3, \ldots $
are independent of each other and independent of $\nu $. But we will admit that functions $F(x) $ and $F^{(d)} (x)$ may depend on the parameter $q$. And let
\begin{eqnarray}\label
 a_1^{(d)} = \mathbf{E}\zeta_1^{(d)} < \infty, \quad  a_1 = \mathbf{E}\zeta_i <
 \infty, \quad  a_2 = \mathbf{E}\zeta_i^2 < \infty, \quad i\geq 2.
  \end{eqnarray}
Such models are of considerable interest for applications (see \cite{Ka97}  and section 4)

Note that, under broad conditions, equality $(\ref{f2.1})$
remains true for geometric sums with delay.

In this article, the results given below for geometric sums with delay will play a central role in the proof of Theorem \ref{t1.1}.

Let's put
\[
 \quad \Delta_{S} (x) =  1 - \exp(-
{x}) - G_{S} \left(\frac{a_1 x}{q^*}\right),   \quad  G_{S} (x) = \mathbf{P}(
S^{(d)}_\nu \leq x).
\]

\begin{lem}\label{l2.1}
Assume that condition (\ref{f2.2}) 
holds for the geometric sum with delay $S^{(d)}_\nu $ and let $x > 0$ be a fixed number. Then, for  $q > 0$ the following estimate holds:
\begin{eqnarray}\label{f2.3}
 \Delta_{S} (x) \geq q^* \exp(- {x}) \left(\frac{a_1^{(d)}}{a_1} -  \frac{a_2}{a_1^2} -g\left(\frac{a_1 x}{q^*}\right)\right)
 ,
   \end{eqnarray}
where
\[
 g(t)= \frac{1}{a_1}\int_t^{\infty} y dF^{(d)}(y) .
\]

\end{lem}

{\it Proof of Lemma \ref{l2.1}.} Let   $N^{(d)} (t)$
denote the counting process constructed from the sequence $(S^{(d)}_k )$, that is,
\[
N^{(d)}(t) = \max(k \geq 1: S^{(d)}_k \leq t),
\]
$H^{(d)} (t)= \mathbf{E}N^{(d)}(t)$ is the renewal function of process
$N^{(d)}(t)$.

Similarly, we introduce a counting process $N(t)$, constructed from a r.v.
$\zeta_k,  k\geq 2$ with a distribution function $F(t)$,  $H (t)= \mathbf{E}N(t)$ is the corresponding renewal function.

Next, we use a simple but important formula
\[
\mathbf{P}(  S^{(d)}_\nu \leq x) = 1- \mathbf{E}(1-q)^{N^{(d)}(x)},
  \]
(see \cite{Ka97}, Ch.3, Theorem 1.2).

We get
\begin{eqnarray}\label{f2.4}
\Delta_{S} (x) =\mathbf{E}(1-q)^{N^{(d)}({a_1 x}/{q^{*}})} - \exp(-
{x}).
  \end{eqnarray}
  
Let's put
\[
 V^{(d)}= V^{(d)}(q, x)= -q^{*}N^{(d)}\left(\frac{a_1 x}{q^{*}}\right) + x
\]
and estimate the value of $ \Delta_{S} (x)$. To do this, we use the well-known elementary inequality
\[
 \forall y \in \mathbb{R}^1  \quad \exp(y) - 1 \geq y.
\]

From this and equality (\ref{f2.4}), we obtain
\begin{eqnarray}\label{f2.5}
 \Delta_{S} (x)  & = & \exp(- {x})(\mathbf{E}  \exp( V^{(d)}) - 1) \geq \exp(- {x}) \mathbf{E}
 V^{(d)} \nonumber \\
 & \geq &  \exp(- {x}) \left(-q^{*}H^{(d)}\left(\frac{a_1 x}{q^{*}}\right) + x\right).
 \end{eqnarray}

Using the law of total probability, we can write the following relation for the function $H^{(d)} (t) $:
\begin{eqnarray}\label{f2.6}
 H^{(d)}(t)  =  \int_0^t (\mathbf{E} N(t-y)+1 ) dF^{(d)}(y) =
 F^{(d)}(t) + \int_0^t H(t-y) dF^{(d)}(y).
  \end{eqnarray}

Next, we present the following uniform estimates for the renewal function  $H(t)$:
\begin{equation} \label{f2.7}
 \frac{t}{a_1} -1 \leq H(t) \leq \frac{t}{a_1}+\frac{a_2}{a_1^2} -1 =
\frac{t}{a_1}+\frac{\sigma^2 }{a_1^2}
 ,
  \end{equation}
(see \cite{lor_70}, \cite{Ka97}, p.57, Proposition 4.2).

From the right-hand inequality in (\ref{f2.7}) and equality (\ref{f2.6}) we have for $t>0$:

\begin{eqnarray}\label{f2.8}
  H^{(d)}(t)  & \leq &  F^{(d)}(t) + \int_0^t \left(\frac{t-y}{a_1}+\frac{a_2}{a_1^2} -1\right) dF^{(d)}(y)  \nonumber \\
& \leq &
  F^{(d)}(t) +\left(\frac{t}{a_1}+\frac{a_2}{a_1^2} -1\right)F^{(d)}(t) - \frac{1}{a_1} \int_0^t y dF^{(d)}(y) \nonumber \\
& \leq &
 \left(\frac{t}{a_1}+\frac{a_2}{a_1^2}\right)F^{(d)}(t) +g(t) - \frac{a^{(d)}_1}{a_1}
 ,
  \end{eqnarray}
where $g(t)=\frac{1}{a_1}\int_t^{\infty} y d F^{(d)}(y)$.

Collecting the estimates (\ref{f2.5}), (\ref{f2.8}), we obtain inequality (\ref{f2.3}) of Lemma \ref{l2.1}.   $ \Box $

\begin{lem}\label{l2.2}
Let condition (\ref{f2.2}) be satisfied for geometric sums with delay $S^{(d)}_\nu $,
$0 < x < 1$ is a fixed number.

Then, for  $q > 0$, the following inequality holds:
\begin{eqnarray}\label{f2.9}
 \Delta_{S} (x) \leq q^* \exp(- {x}) \left(\frac{a_1^{(d)}}{a_1} (1+2x(e-2)) +C_1 (F,x) + C_2 (F,x)\right)
 ,
   \end{eqnarray}
where
\[
 C_1 (F,x) = (2 x^2 (e-2) +x) \left(1- F^{(d)}\left(\frac{a_1 x}{q^{*}}\right)\right),\\
\]
\[
 C_2 (F,x) = (e-2)\left(\left(\frac{2a_2 }{a_1^2}  +1\right) \left(x+\left(\frac{ a_2
  }{a_1^2}-1\right) q^{*}\right)+
  \frac{2a_2 x}{a_1^2 } -2x + q^{*}\right) .
\]
\end{lem}

{\it Proof of Lemma \ref{l2.2}.} To obtain the estimate (\ref{f2.9}) from Lemma \ref{l2.2} we once again use equality (\ref{f2.4}), or rather its corollary: for $0 < x< 1$

\begin{eqnarray}\label{f2.10}
\Delta_{S} (x)  & \leq &  \exp(-x)(\mathbf{E}V^{(d)} +
   (e -2) \mathbf{E}(V^{(d)})^2 )
  \end{eqnarray}
(see \cite{ADM_24}).

Thus, we arrive at the problem of obtaining upper bound estimates for
$\mathbf{E}V^{(d)}$ and $\mathbf{E}(V^{(d)})^2$.

First, we find a lower bound for the renewal function $H^{(d)}
(t)$. To do this, we use the left-hand inequality in (\ref{f2.7}):
\begin{eqnarray}\label{f2.11}
  H^{(d)}(t)  & \geq &  F^{(d)}(t) + \int_0^t \left(\frac{t-x}{a_1} -1\right) dF^{(d)}(x)  \nonumber \\
& \geq &
  F^{(d)}(t) +\left(\frac{t}{a_1} -1\right)F^{(d)}(t) - \frac{1}{a_1} \int_0^t x dF^{(d)}(x) \nonumber \\
& \geq &
 \frac{t}{a_1} - (1-  F^{(d)}(t))\frac{t}{a_1} - \frac{a^{(d)}_1}{a_1}  .
  \end{eqnarray}

The last inequality allows us to provide an upper bound for
$\mathbf{E}V^{(d)}$:
\begin{eqnarray}\label{f2.12}
\mathbf{E}V^{(d)} = -q^{*}H^{(d)}\left(\frac{a_1 x}{q^{*}}\right) + x \leq
q^{*} \frac{a^{(d)}_1}{a_1}+ \left(1- F^{(d)}\left(\frac{a_1 x}{q^{*}}\right)\right){ x}.
  \end{eqnarray}

Next, we consider an upper bound for $\mathbf{E}(V^{(d)})^2$.

Let $H_2^{(d)}(t)= \mathbf{E}|N^{(d)}(t)|^2 $. From the definition, we have

\begin{eqnarray}\label{f2.13}
 \mathbf{E}(V^{(d)})^2 & = & \mathbf{E}\left(-q^{*} N^{(d)}\left(\frac{a_1 x}{q^{*}}\right)+x \right)^2
 =D_1  +D_2,
  \end{eqnarray}
where
\[
D_1 =-2 x q^{*} H^{(d)} \left(\frac{a_1 x}{q^{*}}\right) + 2x^2 , \quad  D_2
=(q^{*})^2 H^{(d)}_2 \left(\frac{a_1 x}{q^{*}}\right) -x^2 .
\]
The estimate for $D_1$ follows directly from inequality (\ref{f2.12}). Indeed,
\begin{eqnarray}\label{f2.14}
  D_1  =
  2x \mathbf{E}V^{(d)}
 \leq 2 x \left(q^{*} \frac{a^{(d)}_1}{a_1}+ \left(1- F^{(d)}\left(\frac{a_1 x}{q^{*}}\right)\right){
 x}\right).
  \end{eqnarray}

To calculate $ D_2 $,  we use the formula (3.10) from \cite{ADM_24}:
\begin{eqnarray}\label{f2.15}
  D_2
    & \leq &  (q^{*})^2 \left(H_2 \left(\frac{a_1 x}{q^{*}}\right)+2 H \left(\frac{a_1
  x}{q^{*}}\right)+1\right) -x^2  \nonumber \\
& \leq &
  (q^{*})^2
  \left(\left( \frac{2\sigma^2 }{a_1^2}  +3\right)H \left(\frac{a_1
  x}{q^{*}}\right)+  \left(\frac{x}{q^{*}}\right)^2
 + \frac{2\sigma^2  x}{a_1^2 q^{*}} + 1\right) - x^2  \nonumber \\
& \leq &
 q^{*}  \left(\left(\frac{2a_2 }{a_1^2}  +1\right) \left(x+\left(\frac{ a_2
  }{a_1^2}-1\right) q^{*}\right)+
  \frac{2a_2 x}{a_1^2 } -2x + q^{*}\right)  .
  \end{eqnarray}

It remains to use estimates (\ref{f2.10}) and (\ref{f2.12}) - (\ref{f2.15}), from which inequality (\ref{f2.9}) of Lemma \ref{l2.2} immediately follows.
 $ \Box $

\section{Proof of Theorem \ref{t1.1}}

Let
\[
\epsilon_{k}(u)= I(\bar{X}_k \geq u) , \quad  I(B) \quad   \mbox{is the indicator of event}\quad   B,
  \]
\[
\nu(u) = \min(k\geq 1 :  \quad \epsilon_{k}(u)= 1).
\]
For cycles of the $2$nd type, we set
\[
 \chi_k =
  \inf(t\geq 0 : \xi_k(t) \geq u /\epsilon_{k}(u)= 1)  ,  \quad  \hat{m}_{\gamma}^{+} = \mathbf{E}(\chi_k^\gamma  /\epsilon_{k}(u)= 1).
\]
That is, $\chi_k $ is the time of the first achievement of level $u$ starting from state $0$ in the $k$th cycle of the  $2$nd type, $\mathbb{\hat{F}}(y)=  \mathbf{P}(\chi_k  < y /\epsilon_{k}(u)= 1)$
is the distribution function of the r.v. $\chi_k $. 

And let
 \begin{eqnarray*}
 \hat{T}_k &= & \left \{
 \begin{array}{rl}
  \chi_k , & \mbox{for} \quad  \epsilon_{k}(u)= 1,  \\
  T_k   , & \mbox{otherwise} .
 \end{array} \right\}
\end{eqnarray*}

The random events
\[
 (  \mathfrak{T}_X (u) \leq x) \quad \mbox{and} \quad \left( \sum_{k=1}^{\nu(u)} \hat{T}_k \leq
x\right)
\]
are equivalent, and therefore
\begin{equation}\label{f3.1}
 G_{X} \left(\frac{x m_1^{-}}{q^*}\right) =
  \mathbf{P}\left( \frac{q^*}{ m_1^{-}} \mathfrak{T}_X (u) \leq x\right)=
\mathbf{P}\left(\frac{q^*}{m_1^{-}}\sum_{k=1}^{\nu(u)} \hat{T}_k \leq x \right).
\end{equation}

The equality of distributions
 \[
 \sum_{k=1}^{\nu(u)-1} \hat{T}_k  \stackrel{def} =
\sum_{k=1}^{\nu(u)-1} {T}^{-}_k
\]
follows directly from the definition and
\[
 \hat{T}_{\nu(u)}  \stackrel{def} = \chi_{\nu(u)} .
\]
Moreover, if the random variables $ \nu(u)$ and $ \hat{T}_k $ are generally dependent, then $ \nu(u)$ and $ {T}^{-}_k , k=1,
2,  \ldots , \nu(u)-1 , $ will be independent. Also, the r.v.
$ \chi_{\nu(u)} $ does not depend on the sum $\sum_{k=1}^{\nu(u)-1} {T}^{-}_k .$

Thus, we have a geometric sum with delay

\[
 \sum_{k=1}^{\nu(u)} \hat{T}_k  \stackrel{def} = \zeta_1^{(d)} +
\sum_{k=2}^{\nu(u)} \zeta_k,
\]
where $ \zeta_1^{(d)} = \chi_{\nu(u)}$,    $\zeta_k={T}^{-}_k , k\geq 2.$

Thus, to prove Theorem \ref{t1.1} we can use Lemmas
\ref{l2.1} and \ref{l2.2}. Let us start with the left inequality in (\ref{f1.3}).

Using the standard properties of conditional expectations, we obtain the following equalities
\begin{eqnarray}\label{f3.2}
m_r =
  (1-q ){m_r}^{-} +q {m_r}^{+},  \quad  \mbox{for}  \quad r=1, 2,
  \gamma.
 \end{eqnarray}

The following estimates follow directly from equalities (\ref{f3.2}): 
 \begin{eqnarray}\label{f3.3}
\frac{1}{1-q }(m_1 - (m_{\gamma })^{1/\gamma}  {q}^{1-1/\gamma})
\leq {m_1}^{-} \leq \frac{m_1}{1-q } ,
\quad m^{+}_1  \leq \left(\frac{m_\gamma}{q}\right)^{1/\gamma} , \nonumber  \\
\frac{1}{1-q }(m_2 -({ m_\gamma })^{2/\gamma}{q}^{1-2/\gamma}) \leq
{m_2}^{-} \leq \frac{m_2}{1-q } ,
\quad m^{+}_2  \leq \left(\frac{m_\gamma}{q}\right)^{2/\gamma} ,  \nonumber \\
\hat{m}_{\gamma}^{+}   \leq
 m_\gamma^{+}\leq \frac{m_\gamma}{q},  \quad 1-\mathbb{\hat{F}}(t)
 \leq \frac{m_\gamma^{+}}{t^\gamma} \leq \frac{m_\gamma}{q t^\gamma}.
\end{eqnarray}

Indeed, for $1<\gamma, r <\infty $ and $1/\gamma +
1/r =1$ from the Holder inequality we have
\begin{eqnarray}\label{f3.33}
m_1^{+}= \mathbf{{E}}(T \mid A_u) = \frac{\mathbf{{E}}(T
\,I(A_u))}{q} \leq \frac{\mathbf{{E}}(T^{\gamma})^{1/\gamma}
q^{1/r}}{q} =\left(\frac{m_\gamma}{q}\right)^{1/\gamma},
\end{eqnarray}
\[
m_1^{-}\geq \frac{m_1 - m_1^{+}q}{1-q} \geq \frac{m_1 -
(\frac{m_{\gamma}}{q})^{1/\gamma}q}{1-q}.
\]
Thus, the first and third inequalities from the first line of (\ref{f3.3})
are established. Analogous estimates can be obtained similarly.

Next, we set in Lemma \ref{l2.1}

\begin{eqnarray}\label{f3.4}
 a_1^{(d)} = \hat{m}_1^{+},
\quad a_1 ={m}_1^{-},  \quad a_2 ={m}_2^{-}, \quad F^{(d)} (x)=
\mathbb{\hat{F}}(x) .
\end{eqnarray}

To derive the left inequality in (\ref{f1.3}) from Lemma \ref{l2.1}, it remains to estimate the function $ g(t)=
\frac{1}{m_1^{-}}\int_t^{\infty} y \, d\mathbb{\hat{F}}(y).$

The required upper bound follows from inequalities (\ref{f3.3}):
\begin{eqnarray}\label{f3.5}
g(t) & =& \frac{1}{{m_1}^{-}}\left( -y(1- \mathbb{\hat{F}}(y))_t^{\infty}
+ \int_t^{\infty} (1- \mathbb{\hat{F}}(y) dy \right) \nonumber \\
& \leq &
  \frac{m_\gamma}{q {m_1}^{-}} \left( \frac{t}{t^{\gamma}} + \int_t^{\infty} \frac{1}{y^{\gamma}}
  dy \right) \nonumber \\
& =&  \frac{\gamma m_\gamma}{q {m_1}^{-} (\gamma -1)t^{\gamma -1}}.
\end{eqnarray}
Substituting $t= x {m_1}^{-}/q^* $ into the last formula from (\ref{f3.4}), we obtain
\[
g\left(\frac{x {m_1}^{-}}{q^*}\right)  \leq \frac{\gamma m_\gamma
(q^{*})^{\gamma - 2}(1+q)}{ ({m_1}^{-})^{\gamma} (\gamma
-1)x^{\gamma -1}}.
\]
This inequality, together with estimate (\ref{f2.3}), gives the left-hand inequality in (\ref{f1.3}).

The proof of the right-hand inequality in (\ref{f1.3})
is primarily based on Lemma  \ref{l2.2}. As above, we also use equalities (\ref{f3.4}) and the estimate

\[
1-\mathbb{\hat{F}}\left(\frac{x {m_1}^{-}}{q^*}\right)  \leq \frac{ m_\gamma
(q^{*})^{\gamma - 1}(1+q)}{ (x{m_1}^{-})^{\gamma} },
\]
which simply follows from inequalities (\ref{f3.3}).

From here and Lema  \ref{l2.2}, we immediately obtain the right-hand inequality in
(\ref{f1.3}).   $ \Box $

Corollary \ref{n1.1} follows directly from Theorem \ref{t1.1}.
Under the conditions of Corollary \ref{n1.1}, the values  $ m_1$, $m_2$, $m_{\gamma}$ are fixed and independent of the parameter
$ u $. Therefore, according to the two-sided estimates for $m^{-}_1$, $m^{-}_2$ in inequalities (\ref{f3.3}), we have:
\[
\quad q(u) \rightarrow 0, \quad  m^{-}_1\rightarrow m_1 , \quad
m^{-}_2 \rightarrow m_2 \quad  \mbox{for}
 \quad u \rightarrow \infty .
\]
Thus, the following implication holds: (\ref{f1.3}) $\Rightarrow$
(\ref{f1.4}). $ \Box $

\begin{cor}\label{z4}
Note that the condition $u \rightarrow\infty$, established in Corollary \ref{n1.1},   is not always satisfied in practice. For example, in the problem of the time to failure of a backup system with recovery from reliability theory, the level $u$ is considered fixed. It is also assumed that $q=q(u, \rho) \rightarrow 0 $ for $ \rho
\rightarrow 0$, where $ \rho$ is the system "load." This means that for highly reliable systems, the average repair time of an element is considered significantly shorter than its average continuous operation time. Clearly, the estimate (\ref{f1.3}) of Theorem \ref{t1.1} can be used in this case as well.
\end{cor}

\section{Examples}

{1. $M/G/1$ queueing system.}

\quad

Consider a single-channel queueing system (QS) receiving a Poisson flow of customers with intensity   $\lambda  $, where $\tau $
denotes the interarrival time, which has an exponential distribution
$\mathbf{P}(\tau < x) = 1- \exp(-\lambda x)$. Let the service time  $\eta $  have an arbitrary distribution
$\mathbf{P}(\eta < x) =G_{\eta} (x) =G(x)$. Suppose that $ \mathbf{E}\eta^k
= b_k < \infty , k\geq 1 $. Thus, in this example, a classical 
${M/G/1}$ queueing system will be studied (see \cite{gk1},
\cite{rio}).

We impose the following condition on the r.v. $\tau$ and $\eta$:
\begin{equation}\label{f4.1}
\rho = \lambda {b_1} <1.
\end{equation}
This condition ensures the existence of a stationary regime in the
queueing system.

We assume that at time $S_0 =0 $, the 
 $M/G/1$ queueing system is empty. At time $\tau_1 $, the first customer arrives, and the first busy period begins, ending at time $S_1 $.
 Accordingly, $S_k $ idenotes the end of the
 $k$-th busy period, $L_k$ its length, and  $L \stackrel{def} = L_k$.

The busy period is the interval starting from the arrival of a request to an empty system and ending at the first subsequent return of the QS to state $0$.

By queue length we mean the total number of requests being processed or awaiting service. Let $Q(t)$ denote the queue length at time $t $.

At time $t $, the queueing system is in state $k $, if  $Q(t)=k . $

Then $Q(t) $ is a regenerative process with regeneration times $S_0 =0, S_1
, S_2 , \ldots $.

Let us consider the problem of estimating the probability of reaching some high level $u$ by process $Q(t)$. To apply Theorem \ref{t1.1} and Corollary \ref{n1.1}, it is necessary to determine the parameters appearing in estimates (\ref{f1.3}), (\ref{f1.4}).

We start with the duration of the $k$-th regeneration period  $T_k $.
It is well known that
 \begin{eqnarray}\label{f4.2}
   m_1 & = &  \mathbf{E}T_k=\mathbf{E}(S_k -S_{k-1})= 1/(1-\rho)\lambda ,  \nonumber \\
 m_2 & = &  \mathbf{E}T_k^2 = \frac{2}{\lambda^2 (1-\rho)} +
 \frac{b_2}{(1-\rho)^3} .
  \end{eqnarray}
   (\cite{kar1}, ch.14, \S 8, \cite{rio}, ch.4, \S 8).

Next we find the probability of the process ${Q}(t)$ exceeding the
level $u$ in one regeneration cycle:
\begin{eqnarray*}
 \quad q = q(u) = \mathbf{P}(\bar{Q}(T_1) \geq u)
  \end{eqnarray*}
  In this section $u$  takes positive integer values.

Let
\[
d_k = \int_{0}^{\infty} \frac{(\lambda x)^{k}}{k!} \exp(-\lambda x)
dG(x)   ,
\]
\[
D_k = \sum_{i=k+1}^{\infty}  d_i   \quad
\]
be the probability that during the service time of one customer, the
QS will receive $ k $ and more than  $ k $ customers, respectively.

An embedded Markov chain $(Q_n), n\geq 0,$ for $M/G/1$ queue is
usually understood as a sequence of values of process $Q(t)$ at service completion times: $Q_n =Q(t_n+0)$, where $t_n$,  $n\geq 1$, denotes the $n$-th service completion time.

Let $ {q}_{k,u} =\mathbf{P}(A_{k,u})$ be a random event
$A_{k,u}$ = \{the embedded Markov chain, starting from the state
$k$ ,  will reach the state  $u$  earlier than the state $0$\} .
More precisely, the values of  $ {q}_{k,u} $ can be determined using Chun's notation for transition probabilities with prohibition (see \cite{chun}, Part 1, $\S 9$):\\
\[
{q}_{k,u} = {}_0q_{k,u}^* =\sum_{n\geq 1} {}_0q_{k,u}^n,  \quad
{}_0q_{k,u}^n =\mathbf{P}(Q_n \geq u, Q_l < u, Q_l \neq 0, l=1, 2,
\ldots, n-1 | Q_0 =k ).
\]
where $ {}_0q_{k,u}^n $ is the probability that chain  $(Q_n) $, starting from state $ k $,  reaches the region $[u, \infty)$ for the first time at step $ n$  without visiting state $ 0$ even once.

The value of $q(u)$ can be determined by solving the following system of linear equations:
\begin{eqnarray}\label{f4.7}
{q}_{u-2,u} & = &  D_1 + d_0 {q}_{u-3,u} + d_1 {q}_{u-2,u}  , \nonumber \\
\ldots & \ldots & \ldots  \quad \ldots  \quad \ldots  \quad , \nonumber \\
{q}_{k,u} & = &  D_{u-k-1} + \sum_{i=0}^{u-k-1} d_i {q}_{i+k-1,u},
\quad 1< k <u-2,
  \nonumber \\
\ldots & \ldots & \ldots  \quad \ldots  \quad \ldots  \quad , \nonumber \\
{q}_{1,u} & = &  D_{u-2} + \sum_{i=1}^{u-2} d_i {q}_{i,u}
 .
  \end{eqnarray}
That is (see \cite{OK_IK_20})
\begin{eqnarray}\label{f4.8}
 \quad q(u) = {q}_{1,u} .
 \end{eqnarray}

Next, we impose additional conditions on the r.v. $\eta$:

\begin{eqnarray}\label{f4.3}
\exists s_0 , 0< s_0 \leq \infty, \quad  \mbox{  such that for any }
\quad
 0< s < s_0 \quad \mathbf{E}\exp(s \eta ) <\infty,
  \end{eqnarray}
\begin{eqnarray}\label{f4.4}
  \mbox{as} \quad
  s \uparrow s_0 \quad \mathbf{E}\exp(s \eta ) \uparrow
  \infty .
  \end{eqnarray}

If the  r.v. $\eta$ has light tails, one can find the asymptotic behavior of the function $q(u)$ for  $u \rightarrow \infty .$
That is, if conditions (\ref{f4.3}), (\ref{f4.4}) hold, then the
following asymptotic formula holds:
\begin{eqnarray}\label{f4.9}
 \quad q(u) = (C_1+o(1)) \exp(-\gamma u) , \quad \gamma =
 \ln(1+{\beta}/{\lambda}),
 \end{eqnarray}
 where   $0< C_1 < \infty$ is a constant independent of $
 u$, and  $\beta >0 $ is the root of the equation
  \begin{eqnarray}\label{f4.10}
\mathbf{E}\exp(\beta  \eta) = 1+\frac{\beta}{\lambda},
  \end{eqnarray}
  (see \cite{ADM_24a}, a similar and well-known result in probability theory can be found in \cite{fe2}, Chap. 12, $\S5$,   \cite{asm}).

The next step is to determine the value of $\hat{m}_1^{+}$ for this example.
As in the general case,  $m_1^{+} $ and $\hat{m}_1^{+}$ denote, respectively, the average duration of a regeneration cycle and the expected time to reach level $u$  by the process  ${Q}(t)$ during second-type regeneration cycles.

As before, we consider transitions of the process $Q(t) $ at service completion times; $(Q_n), n\geq 0$ is the corresponding embedded Markov chain.

Let us introduce the following notation:
\[
\chi_{i,u} =\min(n\geq 1: Q_n \geq u, |Q_0 = i ), \quad
T_{i,u}=\sum_{k =1}^{\chi_{i,u}}\eta_k, \quad {}_0T_{i,u} =T_{i,u}
I(A_{i,u}).
\]
The value  $ {}_0T_{i,u}$ is the time it takes for process  $Q(t) $, starting from state $ i $, to reach the region $[u, \infty)$ provided it never visits state $ 0$; otherwise it equals $ 0 $. It should be noted that the random variable
 \quad $ {}_0T_{i,u}$ defined in this way is either equal to  $0$ or coincides with one of the service end moments. And at the same time it slightly exceeds the actual transition time: $i \rightarrow u $, but not more than $\eta_{\chi_{i,u}} $.

Let's set $  {}_0m_{i,u} =  \mathbf{E} \, {}_0T_{i,u}$, is the number of requests that arrived in the queueing system during the first request $\eta_1$. Next, we write down simple equalities
\begin{eqnarray}\label{f4.01}
\mathbf{E}\eta_1 I(\nu_1 =k)= \int_0^\infty x \frac{(\lambda
x)^k}{k!}\exp(-\lambda x) dG(x) = \frac{k+1}{\lambda} d_{k+1}
 \end{eqnarray}
and
\begin{eqnarray*}
{}_0T_{1,u} =  I(  \nu_1 \geq u-1  )\eta_1 + \sum_{k=1}^{u-2} I(
\nu_1 =k )(\eta_1 +  T_{k,u})I(A_{k,u}).
 \end{eqnarray*}
If we move on to mathematical expectations in the last equality and calculate (\ref{f4.01}), we get

 \begin{eqnarray}\label{f4.02}
{}_0 m_{1,u}  = \sum_{k\geq u-1}\frac{(k+1) d_{k+1}}{\lambda} +
\sum_{k=1}^{u-2}\left(\frac{(k+1)d_{k+1}}{\lambda}q_{k,u}  + d_{k} \,
{}_0 m_{k,u}\right).
 \end{eqnarray}
Describing the possible transitions of the process  $ Q(t)$  after the completion of a service, assuming that initially the queue contains
$i$ requests, $i=2, 3, \ldots , u-2$ we obtain the following system of linear equations:
\begin{eqnarray}\label{f4.03}
{}_0 m_{i,u} & = &  \sum_{k\geq u-i}\frac{(k+1)d_{k+1}}{\lambda} +
 \sum_{k=0}^{u-i-1}\left(\frac{(k+1) d_{k+1}}{\lambda}q_{k+i-1,u} +
d_{k} \, {}_0 m_{i+k -1,u}\right),   \nonumber \\
&  & \quad 1<i < u-2,   \nonumber \\
\ldots & \ldots & \ldots  \quad \ldots  \quad \ldots  \quad  \ldots  \quad  \ldots  \quad  \ldots  \quad ,  \nonumber \\
{}_0  m_{u-2,u} & = & \sum_{k\geq 2}\frac{(k+1)d_{k+1}}{\lambda} +
\frac{ d_{1}}{\lambda}q_{u-3,u}  + \frac{2d_{2}}{\lambda}q_{u-2,u} \nonumber \\
 &+&
d_0 \, {}_0 \, m_{u-3,u}  +  d_1 \, {}_0 \, m_{u-2,u}.
  \end{eqnarray}

Then
\begin{eqnarray}\label{f4.04}
\hat{m}_1^{+} = \frac{1}{\lambda} + \frac{ {}_0 \,m_{1,u}}{q(u)} ,
  \end{eqnarray}
Indeed, the last term in (\ref{f4.04}) represents the expected transition time $1 \rightarrow u $, under the condition that the process  $ Q(t)$ reaches level $u$ before hitting state $0$.

In the case of light-tailed service times, one can obtain a simple upper bound for $\hat{m}_1^{+}$. Consider the functions
\[
M_G (s) = \int_{0}^{\infty} \exp(s x) dG(x) \quad \mbox{ and } \quad
M^{'}_G (s) = \int_{0}^{\infty} x \exp(s x) dG(x).
\]

It is easy to see that if conditions (\ref{f4.1}), (\ref{f4.3})
and (\ref{f4.4}) hold, there exists a number $v_{\lambda}$ such that $ M^{'}_G
(v_{\lambda}) = 1/\lambda .$  Then the duration of the busy period $L$
satisfies the following asymptotic equality:
\begin{eqnarray}\label{f4.5}
\mathbf{P}(L \geq t) \sim C_0 t^{-3/2} \exp(- \alpha t) , \quad
\mbox{ for }t \rightarrow \infty,
\end{eqnarray}
where $C_0 $  is a constant and $\alpha = \lambda +
v_{\lambda} - \lambda M_G (v_{\lambda}) >0$ (see  \cite{cs}, ch.5,
formula (50),  \cite{kyp}).

From the asymptotic relation (\ref{f4.5})  we obtain the following estimate:
 \begin{eqnarray}\label{f4.6}
\hat{m}_1^{+} \leq {m}_1^{+} \leq \frac{1}{\lambda}
+\frac{u(1+o(1))}{\alpha e \, q(u)^{1/u}}.
  \end{eqnarray}

Indeed, to derive (\ref{f4.6}) we proceed as in the proof of (\ref{f3.33}). Namely, for $1<p <\infty $
\begin{eqnarray}\label{f4.61}
\hat{m}_1^{+}\leq \frac{1}{\lambda} + \mathbf{{E}}(L \mid A_u) \leq
\frac{1}{\lambda}
 + \frac{\mathbf{{E}}(L \,I(A_u))}{q(u)} \leq
\frac{1}{\lambda} +\left(\frac{\mathbf{{E}} L^p }{q(u)}\right)^{1/p}.
\end{eqnarray}
Choosing $p = u$ in (\ref{f4.61}) and find the upper bound for the value $\mathbf{{E}} L^u$.

First, we note that (\ref{f4.5}) implies the following estimate: \\
$\exists  t^* >1 $ such that when  $t> t^*$
\[
\mathbf{P}(L > t) \leq (C_0 +1) t^{-3/2} \exp(-\alpha t) .
\]
It follows that
\begin{eqnarray*}
\mathbf{{E}}L^u  & = &  \int_0^\infty \mathbf{P}(L^u  >t) dt \leq t^* +(C_0 +1)\int_{t^*}^\infty t^{-3/2u}\exp(-\alpha t^{1/u}) dt  \nonumber \\
 & \leq &  t^* +(C_0 +1)\int_{t^*}^\infty \exp(-\alpha t^{1/u}) dt  .
  \end{eqnarray*}
Making the substitution $t=(y/\alpha )^u$ in the last integral, we rewrite it as
\begin{eqnarray}\label{f4.62}
  \mathbf{{E}}L^u  & \leq &  t^* +(C_0 +1) \frac{u}{\alpha^u}\int_{y^*}^\infty y^{u-1} \exp(-y) dy  \nonumber \\
 & \leq &  t^* +(C_0 +1) \frac{u }{y^* \alpha^u}\int_{y^*}^\infty y^{u} \exp(-y) dy,
  \end{eqnarray}
where $y^* = (1/\alpha )(t^*)^{1/u} = (1/\alpha )(1+o(1))$.

Next, we will use well-known equalities:
\[
\int_{0}^\infty y^{n} \exp(-y) dy = n! \sim \left(\frac{n}{2\pi}\right)^{1/2}
\left(\frac{n}{e}\right)^n.
\]
Combining these with (\ref{f4.62}), we obtain

\begin{eqnarray}\label{f4.63}
 ( \mathbf{{E}}L^u )^{1/u} & \leq &
  \left( t^* +(C_0 +1) \frac{u }{y^* \alpha^u} u! \right)^{1/u} \sim \frac{u}{e
  \alpha}(1+o(1)).
  \end{eqnarray}
Estimate (\ref{f4.6}) follows directly from (\ref{f4.61}) and
(\ref{f4.63}).

Using (\ref{f4.9}), we rewrite (\ref{f4.6}) in the form:
\begin{eqnarray}\label{f4.11}
\hat{m}_1^{+} \leq    \frac{1}{\lambda}
 +\frac{u \exp(\gamma -1)}{\alpha}(1+ o(1)) .
  \end{eqnarray}

\begin{rem}\label{z5}
(i) Estimate (\ref{f4.11}) slightly sharpens the result from 
\cite{ADM_24a}. Moreover, the proof given in \cite{ADM_24a} appears to contain certain gaps.

(ii) It should be noted that in the case of $M/M/1$ QS   the value of
$\hat{m}_1^{+}$ can be explicitly calculated (see example 2 below). Moreover, it also grows linearly in $u$. Thus, it is not possible to significantly improve estimate (\ref{f4.11}).
\end{rem}

Finally, we determine the mean duration of a type 1 regeneration cycle
$m_1^{-}$ for our example.    Let's put
\[
T_{i,0} =  \min(t>0: \quad Q(t) =0,  |Q(0)
  =i), \quad  {}_uT_{i,0} = T_{i,0} I(  Q(s) < u, 0<s< T_{i,0}
  ). \quad
\]

That is,  ${}_uT_{i,0}$ is the time it takes the process $ Q(t)$ to transition from state $i$ to state  $0$, if it does not reach level $u$, and is equal to 0 otherwise (we assume that at time $t=0$ there are
$i$ requests in the queueing system, one of which has started to be serviced). And let $
 {}_um_{i,0} =\mathbf{E}\, {}_uT_{i,0}.
$

Similarly to (\ref{f4.02}), (\ref{f4.03}), which describe the possible transitions of the process $ Q(t)$ after the end of servicing some request, when initially there are $i$ requests in the QS, $i=1, 2,
\ldots , u-2$ we obtain a system of linear equations for the values  ${}_u
m_{i,0}$:

\begin{eqnarray}\label{f4.14}
{}_u m_{1,0}  & = & \frac{ d_{1}}{\lambda} +
\sum_{k=1}^{u-2}\left(\frac{k+1}{\lambda} d_{k+1} + d_{k} \, {}_u
m_{k,0}\right),  \nonumber \\
\ldots & \ldots & \ldots  \quad \ldots  \quad \ldots  \quad , \nonumber \\
 {}_u m_{i,0} & = & \sum_{k=0}^{u-i-1}\left(\frac{k+1}{\lambda}
d_{k+1} + d_{k} \, {}_u
m_{i+k -1,0}\right), \quad 1<i < u-2,   \nonumber \\
\ldots & \ldots & \ldots  \quad \ldots  \quad \ldots  \quad , \nonumber \\
{}_u  m_{u-2,0} & = & \frac{ d_{1}}{\lambda}  +
\frac{2d_{2}}{\lambda}   +  d_0 \, {}_u \, m_{u-3,0}  +  d_1 \, {}_u
\, m_{u-2,0}.
  \end{eqnarray}

Then
\begin{eqnarray}\label{f4.15}
{m}_1^{-} = \frac{1}{\lambda} + \frac{ {}_u \,m_{1,0}}{1-q(u)} ,
  \end{eqnarray}
where $\frac{1}{\lambda}= E\tau_1 $ is the mean time until the first request arrives, and the last term in the equality (\ref{f4.15}) is in fact the mean duration of the busy period, provided that the process $ Q(t)$ does not reach level $u$ .

Based on Theorem \ref{t1.1} and Corollary \ref{n1.1},
we can formulate our results for the  ${M/G/1}$ QS.

Let us introduce a r.v. $\mathfrak{T}_Q (u) $ as the first time the process $Q (t)$ reaches level $u$, i.e.
$$\mathfrak{T}_Q (u) = \inf(t\geq 0: \quad Q (t) \geq u), $$
and let
\[
 \quad \Delta_Q (x) =  1 - \exp(-
{x}) -   G_{Q} \left(\frac{x m_1^{-}}{q^*}\right),  \quad G_{Q} (x)  =
  \mathbf{P}(  \mathfrak{T}_Q (u) \leq x).
\]

\begin{Stat}\label{tv4.1}
(i) Suppose that condition  (\ref{f4.1}) is satisfied for the
${M/G/1}$ QS, parameter $\lambda$  and function $G(x)$ are fixed,
$0<x< 1$, $u$ takes integer values  and $E\eta^3 = b_3 < \infty$ .

Then, for $u \rightarrow \infty$ the inequalities (\ref{f1.4}) from Corollary \ref{n1.1} hold for $\Delta_Q (x)$,
where  $q=q(u)$ is chosen according to (\ref{f4.8}),
(\ref{f4.9}),  $  m_1, m_2$ are defined by equalities
(\ref{f4.2}), and  $\hat{m}_1^{+}$,  ${m}_1^{-}$ are defined by equalities
(\ref{f4.04}), (\ref{f4.15}), respectively.

(ii) Moreover, if the service time $\eta $ satisfies conditions
(\ref{f4.3}),(\ref{f4.4}), then  $q$ can be chosen according to the formula
(\ref{f4.9}), the value  $\hat{m}_1^{+}$ on the right-hand side of the inequality (\ref{f1.4}) can be replaced by the estimate from
(\ref{f4.11}), where  $ \gamma$ and $ \alpha$  are defined in
(\ref{f4.9}), (\ref{f4.5}), and  ${m}_1^{-} = {m}_1 (1+o(1))$  .

\end{Stat}

{Example 2.  $M/M/1$ queueing system.}

Consider a single-server queueing system that receives a Poisson input (arrival) process with intensity $\lambda
$, and the service time  $\eta $ has an exponential distribution
$$
G(x) = \mathbf{P}(\eta < x) = 1 - \exp(-\mu x), x\ge0 .
$$
In fact, this is a special case of the previous example. Unfortunately, for $M/G/1$ queueing systems  in the general case, we can find some key parameters such as $q=q(u),
\hat{m}_1^{+}, \quad {m}_1^{-}$ only as solutions to certain systems of linear equations. It is clear that for practical use one would like to have simple formulas.

The $ M/M/1$ queueing system is of particular interest because, unlike the general case, it allows explicit computation to find the main parameters in the formulas of Theorem \ref{t1.1} and Corollary \ref{n1.1}.

Let us assume that for this queueing system the condition
(\ref{f4.1}) is true, and $b_1 =1/\mu $ .  Then the process  $Q(t) $ is regenerative with moments of regeneration $S_0 , S_1, \ldots $ (see example 1). It remains to find the parameters in the formulas of Theorem \ref{t1.1} and Corollary \ref{n1.1}.

Our problem is closely related to the classical problem of a random walk with two-sided absorption $x=0$ and $x=u$, and probabilities of up $0< p=\lambda/(\lambda +\mu)< 1/2$ and down $1-p=\mu/(\lambda +\mu)$ transitions (\cite{fe1}, ch.14). Namely, let us consider the embedded Markov chain for queueing systems  $Q^*_n =
Q(t^*_n) $, where $t^*_n $ are the moments of change of states  $ Q(t) $.

Then, if we choose the initial state $x_0 =1$, then the probability that the Markov chain $Q^*_n$  at the  $n$th step will reach the point $x=0$ and will not reach the level $u$, is equal to the probability of absorption at the point $x=0$ of such a random walk, which also starts from the same point $x_0=1$. A similar situation will occur when we consider the absorption probabilities at point $x=u$.

Therefore, we will first present some results for such a random walk. Let \, ${}_u f_{1,0}^{(n)}$ \, and
\,
  ${}_0 f_{1,u}^{(n)}$\, denote the probabilities that a bounded random walk, starting from point $x_0=1$,  is absorbed at step $n$th at point $x=0$ and  $x=u$, respectively. And let
 \[
 P_{1,0}(s)= \sum_{n=0}^{\infty} {}_u f_{1,0}^{(n)} s^n       ,
\quad  P_{1,u}(s)= \sum_{n=0}^{\infty}  {}_0 f_{1,u}^{(n)} s^n \quad
-
 \]
be the generating functions of the distributions  $({}_u f_{1,0}^{(n)})$ and 
 $({}_0 f_{1,u}^{(n)})$.

It is known  (\cite{fe1}, ch.14, \S 4) that under the condition
$2\sqrt{p(1-p)} \, s < 1$
\begin{equation}\label{f4.22}
P_{1,0}(s)= \frac{1-p}{p} \, \frac{h_1^{u-1}(s)-
h_2^{u-1}(s)}{h_1^{u}(s)- h_2^{u}(s)} ,
\end{equation}

\begin{equation}\label{f4.23}
P_{1,u}(s)=  \frac{h_1 (s)- h_2 (s)}{h_1^{u}(s)- h_2^{u}(s)} ,
\end{equation}
where
\[
h_1 (s)= \frac{1+(1-4p(1-p)s^2)^{1/2}}{2ps} ,  \quad h_2 (s)=
\frac{1 -(1-4p(1-p)s^2)^{1/2}}{2ps} .
\]

Next, for the random walk we compute the probability of absorption at the upper boundary $x=u$ when it starts from $x_0=1$. We have
\begin{equation}\label{f4.24}
h_1 (1)= \frac{1-p}{p}= \frac{1}{\rho} ,  \quad h_2 (1)= 1 , \quad
\rho =\frac{ \lambda}{\mu},
\end{equation}
\[
P_{1,u}(1)=  \frac{h_1 (1)- h_2 (1)}{h_1^{u}(1)- h_2^{u}(1)} =
\frac{1/\rho  -1}{1/\rho^u -1} =\frac{(1-\rho) \rho^{u-1}}{1 -\rho^u
},
\]
which usually coincides with known results (\cite{fe1}, ch.14,
\S 2, the case of the birth and death process see \cite{ADM_24}).
It is clear that then for the $ M/M/1$ QS
\begin{equation}\label{f4.25}
q=q(u)= P_{1,u}(1)= \frac{(1-\rho) \rho^{u-1}}{1 -\rho^u }.
\end{equation}

Next, we consider the expected time for the random walk to reach absorption at the upper boundary $x=u$ starting from $x_0=1$. Let  $\varepsilon $ be a random variable taking values
$+1$ and $-1$ with probabilities $p$ and $1-p$ respectively. Let
$(\varepsilon_k)$ be independent copies of $\varepsilon $, $s_n =1 +
\sum_{k=1}^n \varepsilon_k$,

\begin{eqnarray*}
{}_0 \chi_{1,u} &= & \left \{
 \begin{array}{rl}
\min\{n\geq1: \quad  s_n \geq u, \quad s_k >0, k=1, 2,\ldots, n-1\},  & \mbox{if such n exists,}     \\
 0  , & \mbox{if no such n exists}.
 \end{array} \right\}
\end{eqnarray*}

That is, we want to compute the expected value:
\begin{eqnarray}\label{f4.26}
{}_0 m^*_{1,u} & = & \mathbf{E} \,{}_0 \chi_{1,u}  =
\sum_{n=0}^{\infty}n \, {}_0 f_{1,u}^{(n)} =
P^{'}_{1,u}(1)=  \nonumber \\
 & = &
\frac{(h^{'}_1 (1)- h^{'}_2 (1))(h^{u}_1 (1)- h^{u}_2 (1))
-u(h^{u-1}_1 (1)h^{'}_1 (1)- h^{u-1}_2 (1)h^{'}_2 (1) )}
 {(h_1^{u}(1)- h_2^{u}(1))^2}
 .
  \end{eqnarray}

Elementary calculations allow us to find the values that are included in the last equality:
\begin{equation}\label{f4.27}
h^{'}_1 (1)- h^{'}_2 (1)=
\left(\frac{(1-4p(1-p)s^2)^{1/2}}{ps}\right)^{'}_{s=1} =\frac{ -4(1-p)}{1-2p} -
\frac{1-2p}{p} = \frac{-1}{p(1-2p)}.
\end{equation}
Similarly
\begin{equation}\label{f4.28}
h^{'}_1 (1)=  - \frac{1-p}{p(1-2p)}, \quad h^{'}_2 (1)=
\frac{1}{1-2p} .
\end{equation}
Putting together equalities (\ref{f4.24}), (\ref{f4.26})-(\ref{f4.28})
we get
\begin{eqnarray}\label{f4.29}
{}_0 m^*_{1,u} & = &  \nonumber \\
 & = &
\frac{(-1/p(1-2p))(1/\rho^{u}- 1) -u((p-1)/(\rho^{u-1} p(1-2p)) -
1/(1-2p) )}
 {(1/\rho^{u}- 1)^2}  \nonumber \\
 & = &
\frac{ u(1/\rho^{u}+1) -(1/p)(1/\rho^{u}- 1)}
 {(1-2p)(1/\rho^{u}- 1)^2}
  .
  \end{eqnarray}

Next, we transition from the expected value ${}_0 m^*_{1,u}$ of the random walk to the expected value ${}_0 m_{1,u}$ of the process $Q(t)$.

Let us denote by $ \tau_{i}^* $  the duration of stay of process
$Q(t)$ in state $i$. It is well known that $ \tau_{i}^* $ has an exponential distribution and $\mathbf{E} \tau_{i}^* =1/(\lambda +
\mu), \quad i\geq 1 $.

Let, as above, ${}_0 T_{1,u}$ be the time of transition of process
$Q(t)$ from state $1$ to state $u$ without reaching $0$. Then
\[
{}_0 T_{1,u}= \sum_{k=1}^{{}_0 \chi_{1,u}}\tau_{i_k}^* ,
\]
where ${i_k \neq 0}$ runs through the corresponding stages through which the transition is made: $1 \rightarrow u.$

According to the definition, the r.v. ${}_0 \chi_{1,u}$  s completely determined by the sequence $(\varepsilon_k )$, which does not depend on the r.v.
$(\tau_{i}^* )$. Therefore, by the Wald identity and using estimate
(\ref{f4.29}), we obtain
\begin{eqnarray}\label{f4.30}
{}_0 m_{1,u} & = & \mathbf{E} \,{}_0 T_{1,u}=\mathbf{E}
\sum_{k=1}^{{}_0 \chi_{1,u}}\tau_{i_k}^* ={}_0 m^*_{1,u}  \frac{1}{\lambda + \mu}\nonumber \\
 & = &
\frac{ u(1/\rho^{u}+1) -(1/p)(1/\rho^{u}- 1)}
 {(1-2p)(1/\rho^{u}- 1)^2} \frac{1}{\lambda +
\mu}
 .
  \end{eqnarray}

Using equalities (\ref{f4.25}) and (\ref{f4.30}) we can now calculate parameter  $\hat{m}_1^{+}$, which appears in inequality (\ref{f1.4}):
\begin{eqnarray}\label{f4.31}
\hat{m}_1^{+} & = & \frac{1}{\lambda} +\frac{ {}_0 m_{1,u}}{q(u)} \nonumber \\
 & = &  \frac{1}{\lambda} + \frac{ u(1/\rho^{u}+1) -(1/p)(1/\rho^{u}- 1)}
 {
(1-2p)(1/\rho^{u}- 1)^2} \frac{\rho}{(\lambda +
\mu)(1-\rho)\rho^u} \nonumber \\
 & = &
 \frac{1}{\lambda} + \frac{\rho( u  -1/p +O(u\rho^{u}))}
 {
(1-2p)(1 - \rho )(\lambda + \mu)}
   .
  \end{eqnarray}

Next, for process $Q(t)$, we find the expected duration of a type-$1$ regeneration cycle, ${m}_1^{-}$. Introduce the following notations:
$ {}_0 \,m^{*}_{1,u} = P^{'}_{1,u}(1)$,    ${}_u
\,m^{*}_{1,0} =P^{'}_{1,0}(1)$,
 $m^*$ be the expected absorption time in states $x=0$ or $x=u$  for a bounded random walk that starts from point $1$. It is well known (see \cite{fe1},
ch.14, \S 3) that
\begin{eqnarray}\label{f4.32}
{}_0 \,m^{*}_{1,u} + {}_u \,m^{*}_{1,0} =  m^* = \frac{1}{1-2p}
\left(1-\frac{u(2-1/p)}{1-(1/p -1)^u}\right) .
\end{eqnarray}
In addition, as above (see equality (\ref{f4.30})), we have
\begin{eqnarray}\label{f4.33}
 {}_u m_{1,0}  = {}_u m^*_{1,0} \frac{1}{\lambda + \mu}.
\end{eqnarray}

Thus
\begin{eqnarray}\label{f4.34}
{m}_1^{-} = \frac{1}{\lambda} +\frac{ {}_u \,m_{1,0}}{(1-q(u))}=
\frac{ {}_u \,m^{*}_{1,0}}{(\lambda + \mu)(1-q(u))} .
  \end{eqnarray}
Since the value ${}_0 \,m^{*}_{1,u}$ is calculated in (\ref{f4.29}), equalities (\ref{f4.32}) and (\ref{f4.34}) allow us to find the exact value of ${m}_1^{-} $.

Note that for $ M/M/1$ queues  $\E\eta^2 =
2/\mu^2$, and therefore
\begin{eqnarray}\label{f4.35}
\frac{m_2}{m_1^2} =2(1-\rho) +\frac{2\rho^2}{1-\rho} .
 \end{eqnarray}
 
Combining the above estimates with Corollary \ref{n1.1} yield the following statement.

\begin{Stat}\label{tv4.2}
(i) Suppose that condition  (\ref{f4.1}) is satisfied for the
${M/M/1}$ queue, parameters $\lambda$ and $\mu$ are fixed,  $0<x< 1$,
$u$ takes integer values.

 Then for $u \rightarrow \infty$
\begin{eqnarray}\label{f4.36}
   &&\exp(-x)(1-\rho)\rho^{u-1}\left(\frac{\hat{m}_1^{+}}{m^{-}_1} -2(1-\rho +\frac{\rho^2}{1-\rho})+o(1)\right)
    \leq  \Delta_Q (x)  \leq    \exp(-x)(1-\rho)\rho^{u-1}\times   \nonumber \\
   & \times &
   \left(\frac{\hat{m}_1^{+}}{m_1^{-}}(1+2x(e-2))
  +  8(e-2)\left(1-\rho +\frac{\rho^2}{1-\rho}\right)x -x +o(1) \right),
  \end{eqnarray}
where  $\hat{m}_1^{+}$ and  ${m}_1^{-}$ are given by equalities
(\ref{f4.31}) and (\ref{f4.34}) respectively.

\end{Stat}

\end{document}